\magnification=\magstep1
\input amstex
\documentstyle{amsppt}
\hsize = 6.5 truein
\vsize = 9 truein

\UseAMSsymbols

\define \J {\Bbb J}
\define \rL {\Bbb L}
\define \K {\Bbb K}
\define \rO {\Bbb O}
\define \N {\Bbb N}
\define \AN {\Bbb N^*}
\define \rP {\Bbb P}
\define \CP {\Cal P}
\define \CET {\Cal C (ET)}
\define \s {\Bbb S}
\define \zs {\Bbb S_0}
\define \CRET {\Cal C (RET)}
\define \Z {\Bbb Z}
\define \Zi {\Bbb Z[i]}
\define \Q {\Bbb Q}
\define \R {\Bbb R}
\define \C {\Bbb C}
\define \oN {\overline {\Bbb N}}
\define \CPET {\Cal C (PET)}
\define \CePET {\Cal C_e (PET)}
\define \oR {\overline {\Bbb R}}
\define \RP {\R^{+*}}
\define \cal {\text {cal}}
\define \th	{\theta}
\define \cp {\text {cap}}
\define \gri {\text {gri}}
\define \lgr {\text {lgr}}
\define \grs {\text {grs}}
\define \ugr {\text {ugr}}
\define \gr {\text {gr}}
\define \e	{\epsilon}
\define \de	{\delta}
\define \p {\Bbb P}
\define \Gi {\Bbb Z[i]}

\define \ee {\'e}
\define \ea {\`e}
\define \oo {\^o}

\define \CC {\Cal C}
\define \CDD {\Cal D}
\define \CF {\Cal F}
\define \CG {\Cal G}
\define \CH {\Cal H}
\define \CI {\Cal I}
\define \CL {\Cal L}
\define \CN {\Cal N}
\define \CT {\Cal T}
\define \CU {\Cal U}
\define \CV {\Cal V}
\define \CW {\Cal W}
\define \CA {\Cal A}
\define \CB {\Cal B}
\define \CE {\Cal E}
\define \CR {\Cal R}
\define \CS {\Cal S}
\define \CM {\Cal M}
\define \CO {\Cal O}
\define \CX {\Cal X}
\define \CZ {\Cal Z}

\define \card {\text{card}}

\define \bS {\bf S \rm}

\define \noo {$\text{n}^o$}

\define \T- {\overset {-1}\to T}
\define \G- {\overset {-1}\to G}
\define \GA- {\overset {-1}\to \Gamma}
\define \RA- {\overset {-1}\to R}
\define \PA- {\overset {-1}\to P}
\define \ig {g^{-1}}

\define \Ž {\'e}
\define \ {\`e}
\define \ˆ {\`a}
\define \ {\`u}

\define \™ {\^o}
\define \ {\^e}
\define \" {\^\i}
\define \ž {\^u}
\define \‰ {\^a}

\define \rd {\roman d}
\define \rD {\roman D}
\define \rad {\roman {ad}}
\define \rp {\roman p}
\define \rs {\roman s}

\define \su {\subheading}
\define \noi {\noindent}

\UseAMSsymbols
\topmatter
\title Produit d'entrelacement et action triangulaire
d'alg\`ebres de Lie
\endtitle

\author Barben-Jean COFFI-NKETSIA et Labib HADDAD \endauthor

\address {120 rue de Charonne, 75011 Paris, France;
e-mail: labib.haddad\@wanadoo.fr } \endaddress

\endtopmatter
\document
\refstyle{A}
\widestnumber\key{ABCD}

\head R\'esum\'e\endhead

\

En mimant les lois d'op\Ž rations infinit\Ž simales des
alg\ bres de Lie sur les vari\Ž t\Ž s analytiques
banachiques, on introduit de mani\ re purement alg\Ž brique
la notion d'action formelle d'une alg\ bre de Lie sur un
espace vectoriel. Ensuite, par analogie avec le cas des
groupes abstraits, et en faisant op\Ž rer les alg\ bres de
Lie \lq\lq en cascade", on d\Ž finit produit
d'entrelacement (\lq\lq wreath product") et action
triangulaire pour les alg\ bres de Lie. On d\Ž montre enfin
un th\Ž or\ me du type Kaloujnine-Krasner pour les
extensions d'alg\ bres de Lie.

\

\head Introduction\endhead

\

Le produit d'entrelacement (\lq\lq
wreath product") $W$ de deux groupes quelconques se d\Ž finit
commod\Ž ment en faisant agir ces deux groupes \lq\lq en
cascade", ce qui conduit \ˆ \ la notion classique d'action
triangulaire. Lorsque l'on a affaire \ˆ \ des groupes \bf de
Lie\rm, il est naturel de vouloir munir $W$ (ou du moins un
de ses sous-groupes substantiels) d'une sturcture convenable
de groupe de Lie.

\

\bf \lq\lq Peut-on d\Ž finir le produit d'entrelacement de
deux groupes de Lie ?"\rm

\

Cette question, pos\Ž e par M. Krasner au premier de nous
deux, est \ˆ \ l'origine de ce travail.

\

Le probl\ me n'est pas ais\Ž \ et ne semble toujours pas
avoir \Ž t\Ž \ r\Ž solu.

\

Afin de tourner la difficult\Ž , et comme premi\ re \Ž
tape, nous avons voulu chercher d'abord une bonne d\Ž
finition pour un produit d'entrelacement de deux \bf
alg\ bres \rm de Lie. C'est cela que nous proposons
ci-dessous.

\

Pour y parvenir, nous avions besoin d'une notion d'\bf
action \rm pour une alg\ bre de Lie quelconque. Elle devait
\ tre suffisamment g\Ž n\Ž rale pour \ tre utilisable dans
notre contexte. Nous avons \Ž t\Ž , ainsi, amen\Ž s \ˆ \
forger un certain nombre d'outils. Ils nous paraissent
avoir, par ailleurs, de l'int\Ž r\ t en eux-m\ mes. Nous en
pr\Ž sentons l'essentiel dans cette note.

\

Nous avons ainsi d\Ž fini, notamment, la d\Ž rivation
suivant une s\Ž rie formelle \ˆ \ variables et coefficients
dans un espace vectoriel $X$ quelconque. Ce qui m\ ne \ˆ \
l'introduction d'un objet nouveau, \bf l'alg\ bre de
Lie $S(X)$ des s\Ž ries formelles sur $X$\rm. D'o\ \ la
notion nouvelle d'action formelle d'une alg\ bre de Lie
quelconque sur l'espace vectoriel $X$. Et, singuli\ rement,
l'action fondamentale d'une alg\ bre de Lie sur elle-m\ me
(qui m\Ž rite d'\ tre signal\Ž e). D'o\ \ d\Ž coule alors,
assez naturellement, notre d\Ž finition du produit
d'entrelacement.

\

Un th\Ž or\ me de repr\Ž sentation (\sl \ˆ \ la
Kaloujnine-Krasner\rm, voir [6]) vient enfin, \ˆ \ point
nomm\Ž , illustrer le bon fonctionnement de ce produit
d'entrelacement : toute extension $C$ d'une alg\ bre de Lie
$B$ par une alg\ bre de Lie $A$ se plonge dans le produit
d'entrelacement de l'alg\ bre $B$ par l'alg\ bre $A$. On
remarquera, en particulier, la formule de ce plongement,
(3) au paragraphe \bf 13\rm.

\

Les constructions, ici, sont tr\ s diff\Ž rentes et
beaucoup plus complexes que dans le cas classique des
groupes. Ce pourrait \ tre cependant le premier pas vers une
\lq\lq bonne" d\Ž finition du produit d'entrelacement pour
les
\bf groupes de Lie\rm.

\ 

\{\bf En effet\rm, on peut penser que
l'alg\ bre de Lie $L(W)$ d'un produit d'entrelacement 
\sl convenable \rm $W$ de deux groupes de Lie $G$ et $H$
devrait
\ tre le produit d'entrelacement de leurs alg\ bres de Lie
respectives
$L(G)$ et $L(H)$. L'on pourrait ainsi esp\Ž rer remonter
du produit d'entrelacement des alg\ bres \ˆ \ celui des
groupes de Lie. Une des difficult\Ž s qui se pr\Ž sentent
alors c'est le passage oblig\Ž \ du fini
\ˆ \ l'infini car, m\ me si les deux alg\ bres de Lie
$L(G)$ et
$L(H)$ sont, toutes deux, de dimension finie, leur produit
d'entrelacement a une dimension infinie. Le produit
d'entrelacement des deux groupes $G$ et $H$ devrait alors \
tre \sl model\Ž \rm
\ sur un espace de Hilbert, ou de Banach
pour le moins, de dimension infinie.\}

\

Nous pr\Ž voyons de faire para\^\i tre ult\Ž rieurement
tous les d\Ž tails de nos constructions pour les alg\ bres
de Lie. Ils sont nombreux. Ils ne sont pas toujours imm\Ž
diats. Un texte est en pr\Ž paration. Les num\Ž ros  de la
forme
$<n>$ y renvoient. 

\

Nous pr\Ž sentons, ici, les \sl grandes lignes \rm de la d\Ž
marche.

\

\head D\Ž veloppement \endhead

\

On se fixe un corps commutatif $K$ de caract\Ž ristique
nulle.

\

\bf Tous les espaces vectoriels et toutes les alg\ bres de
Lie consid\Ž r\Ž s sont suppos\Ž s avoir $K$ comme corps des
scalaires.\rm 

\

On d\Ž signe par $E,F,X,Y,$ des espaces vectoriels et par
$A,B,C,$ des alg\ bres de Lie. Par $m,n,r,$ on d\Ž signera
des entiers naturels quelconques.

\

Pour chaque $m$, on d\Ž signera par $L_m(E;F)$ l'ensemble
des applications $m$-lin\Ž aires sur $E$ \ˆ \ valeurs dans
$F$.

\

\su{1 S\Ž ries formelles} On dira qu'une application $f : E
\to F$ est un \bf polyn\™ me homog\ ne de degr\Ž \ $m$ \ˆ \
variables dans $E$ et coefficients dans $F$ \rm lorsqu'il
existe $u \in L_m(E;F)$ tel que
$$f(x) = u(x,x,\dots,x) \ \ \text{pour tout} \ \ x \in E.$$

On dira alors que $f$ est \bf d\Ž termin\Ž \ \rm  par $u$.

\

On d\Ž signera par $F[E]_m$ l'ensemble de ces polyn\™ mes
homog\ nes de degr\Ž \ $m$. C'est naturellement un espace
vectoriel (sur $K$). On posera
$$F[E] = \underset m\to \bigoplus F[E]_m \ , \ F[[E]] =
\prod_m F[E]_m.$$
On appelle alors \bf polyn\™ me \rm (resp. \bf s\Ž rie
formelle\rm) \ˆ \ variables dans $E$ et coefficients dans $F$
tout \Ž l\Ž ment de $F[E]$ (resp. de F[[E]]).

\

\su{Remarque} Dans le cas o\ \ $E$ et $F$ sont des espaces
norm\Ž s, on d\Ž signe par $\hat P(E;F)$ l'ensemble des \bf
s\Ž ries formelles \ˆ \ composantes continues \rm sur $E$ \ˆ
\ valeurs dans $F$ (voir Bourbaki [1], p.88-89). Ainsi $\hat
P(E;F)$ est un sous-espace vectoriel de $F[[E]]$, qui lui
est \Ž gal lorsque la dimension de $E$ est \bf finie\rm.

\

\su{2 Sym\Ž trisation} Pour chaque 
$$u \in L_m(E;F) \ , \ z = (z_1,\dots,z_r) \in E^r \ , \ p =
(p_1,\dots,p_r) \in \Z^r,$$
si $m = p_1 + \dots + p_r$, on d\Ž signe par $\tilde
u(z;p)$ la somme de tous termes de la forme
$u(x_1,\dots,x_m)$ o\ , parmi les $x_i$ $(i = 1,\dots,m)$,
il y en a exactement $p_j$ qui sont \Ž gaux \ˆ \ $z_j$, pour
$j = 1, \dots,r$. S'il n'existe aucun terme de cette
forme, on convient que $\tilde
u(z;p) = 0$.

\

On a alors le r\Ž sultat suivant.

\

\su{3 Th\Ž or\ me} $<1>$. Soient $u$ et $v$ des \Ž l\Ž
ments de $L_m(E;F)$. On suppose que $u(x,\dots,x) =
v(x,\dots,x)$ pour tout $x \in E$ (autrement dit, $u$ et $v$
d\Ž terminent le m\ me polyn\™ me homog\ ne). Alors
$$\tilde u(z;p) = \tilde v(z;p) \ \ \text{pour tous} \ \ z
\in E^r \ , \ p \in \Z^r.$$
On \Ž tablit ce th\Ž or\ me en utilisant le r\Ž sultat
suivant.

\

Pour $p = (p_1,\dots,p_r) \in \N^r$ et $t = (t_1,\dots,t_r)
\in K^r$, on pose $t^p = t_1^{p_1} \dots t_r^{p_r}$,
convention des multiindices, et $|p| = p_1 ± \dots + p_r$. 

\

\su{4 Lemme} $<2>$. Soient $m,r,$ des entiers naturels et,
pour chaque $p \in \N^r$ tel que $|p| \leq m$, soit $a_p$ un
\Ž l\Ž ment de $E$. Si $\dsize\sum_{|p| \leq m} t^pa_p = 0$
pour tout $t \in K^r$, alors $a_p = 0$ pour tout $|p| \leq
m$.

\

\su{5 D\Ž rivation suivant une s\Ž rie formelle}

\

Soient $\xi \in X[X]_r$ et $f \in F[X]_m$.

\  

On suppose que $f$ est \bf d\Ž termin\Ž \ \rm  par $u \in
L_m[X;F)$. Pour chaque $x \in X$, on d\Ž signe par $(\xi
f)(x)$ la somme de tous les termes de la forme
$u(x_1,\dots,x_m)$ o\ , parmi les $x_i$ $(i = 1,\dots,m)$,
un seul est \Ž gal \ˆ \ $\xi(x)$ et les autres sont \Ž gaux
\ˆ \ $x$. Autrement dit, 
$$(\xi f)(x) = \tilde
u((\xi(x),x);(1,m-1)).$$ 
Cela d\Ž finit une application $\xi f : X \to F$.

\

On montre (proposition $<3>$), que $\xi f$ est un polyn\™ me
homog\ ne. Il ne d\Ž pend pas du choix de $u$, d'apr\ s le
th\Ž or\ me 3, et l'on a $\xi f \in F[X]_s$ o\ \ $s = r +
m - 1$, avec la convention $F[X]_s = \{0\}$ pour $s < 0$.

\

\

\su{D\Ž finition} On \Ž crira $S(X)$ au lieu
de $X[[X]]$.

\

Soient $\xi = (\xi_r) \in S(X)$ et $f = (f_m)
\in F[[X]]$ des s\Ž ries formelles. Pour tout $s \geq 0$,
posons
$$g_s = \sum_{r + m - 1 = s} \xi_r f_m.$$
On d\Ž signe par $\xi f$ la s\Ž rie formelle $(g_s) \in
F[[x]]$ et on l'appelle \bf d\Ž riv\Ž e de $f$ suivant
$\xi$\rm.

\
 
\su{6 L'alg\ bre de Lie $S(X)$}  Etant donn\Ž es des s\Ž
ries formelles $\xi$ et $\eta$ dans
$S(X)$, on peut consid\Ž rer $\xi \eta$ la d\Ž riv\Ž e de
$\eta$ suivant $\xi$, et $\eta \xi$ la d\Ž riv\Ž e de $\xi$
suivant $\eta$. On pose $[\xi,\eta] = \xi \eta - \eta \xi$.

\

\bf On montre alors\rm, (th\Ž or\ me $<4>$), \bf que
l'espace vectoriel
$S(X)$ muni \rm du crochet \bf ainsi d\Ž fini est une alg\
bre de Lie\rm.

\

\su{7 Action d'une alg\ bre de Lie sur un espace vectoriel}
On appellera \bf action formelle \rm (\ˆ \ droite) de
l'alg\ bre de Lie $A$ sur l'espace vectoriel $X$ tout
homomorphisme d'alg\ bres de Lie $\rD : A \to S(X)$. Ainsi,
on fait agir l'alg\ bre de Lie $A$ sur l'espace vectoriel
$X$ \sl au travers \rm de son alg\ bre de Lie $S(X)$.

\

\su{Exemple originel} On suppose que $K = \R$ ou $\C$, que
$A$ est une alg\ bre de Lie normable compl\ te et que $X$
est un espace de Banach. On consid\ re un voisinage ouvert
$U$ de $0$ dans $X$, et une \bf loi d'op\Ž ration infinit\Ž
simale \ˆ \ droite\rm, analytique, $a \mapsto \rD_a$, de $A$
dans la vari\Ž t\Ž \ analytique $U$ (voir Bourbaki [4],
p.139).

\

Ainsi, pour chaque $a \in A$, le champ de vecteurs $\rD_a$
sur
$U$ est analytique. Bien entendu, l'injection canonique $h :
U \to X$ est analytique. On consid\ re l'application
$\rD_a(h) : U \to X$ (voir Bourbaki [2], 8.2.2 et 8.2.3,
p.10).

\

Elle est analytique donc repr\Ž sentable au voisinage de
l'origine par une s\Ž rie formelle (convergente) \ˆ \
composantes continues, c'est-\ˆ -dire par un \Ž l\Ž ment de
$\hat P(X;X) \subset S(X)$, que nous d\Ž signerons encore par
$\rD_a$.

\

\bf On v\Ž rifie \rm (th\Ž or\ me $<5>$) \bf que
l'application
$\rD : A
\to S(X)$ ainsi d\Ž finie est une action formelle de $A$
sur $X$\rm. On dira que c'est l'action formelle
\bf d\Ž duite \rm de la loi d'op\Ž ration infinit\Ž simale
donn\Ž e.

\

\su{8 Produit d'entrelacement} On consid\ re une action
formelle $\rD$ de $A$ sur $X$ et une action
formelle $\rd$ de $B$ sur $Y$. On consid\ re les
deux espaces vectoriels produits
$$W = A[[Y]] \times B \ \ \text{et} \ \  \ Z = X \times Y.$$

\

On va d\Ž finir d'abord une structure d'alg\ bre de Lie sur
$W$, que l'on appellera \bf produit d'entrelacement\rm, puis
une action formelle $\Delta : W \to S(Z)$, que
l'on appellera \bf action triangulaire\rm.

\

\su{8.1 Crochet sur $A[[Y]]$} 

\

Pour $f \in A[Y]_n$ et $g \in
A[Y]_r$ et chaque $y \in Y$, on pose $[f,g](y) =
[f(y),g(y)]$. Cela d\Ž finit une application $[f,g] : Y \to
A$. On v\Ž rifie (lemme $<6>$) que $[f,g] \in A[Y]_{n + r}$.
Plus g\Ž n\Ž ralement, pour $f = (f_n) \in A[[Y]]$ et $g =
(g_r) \in A[[Y]]$, on pose $[f,g]_s = \sum_{n+r =s}
[f_n,g_r]$ et enfin $[f,g] =([f,g]_s) \in A[[Y]]$.

\

Le \bf crochet \rm ainsi d\Ž fini munit l'espace vectoriel
$A[[Y]]$ d'une stucture d'alg\ bre de Lie \bf h\Ž rit\Ž e
\rm de celle de $A$ (proposition $<7>$.)

\

D\Ž signons par $\frak d (A[[Y]])$ l'alg\ bre de Lie des
d\Ž rivations de l'alg\ bre de Lie $A[[Y]]$.

\

\su{8.2 Homomorphisme de $B$ dans $\frak d (A[[Y]])$}

\

Reprenons l'action formelle $\rd : B \to S(Y)$.
Pour chaque $b \in B$, on a $\rd_b \in S(Y)$; et pour chaque
$a \in A[[Y]]$, la d\Ž riv\Ž e $\rd_b a$ de $a$ suivant
$\rd_b$ appartient \ˆ \ $A[[Y]]$ (voir ci-dessus, au \bf
5\rm).

\

Ainsi $\rd_b$ d\Ž finit une application de $A[[Y]]$ dans
elle-m\ me. On v\Ž rifie (proposition $<8>$) que cette
application est une d\Ž rivation de l'alg\ bre de Lie
$A[[Y]]$ et on la d\Ž signe par $\sigma(b)$.

\

On obtient ainsi une application $\sigma : B \to \frak d
(A[[Y]])$ et on v\Ž rifie (proposition $<9>$) que $\sigma$
est un homomorphisme d'alg\ bres de Lie.

\

\su{8.3 Crochet sur $W = A[[Y]] \times B$}

\

Soient $(a,b)$ et $(a',b')$ des \Ž l\Ž ments de $W$.

\

On pose
$$[(a,b),(a',b')] = ([a,a'] + \rd_b a' - \rd_{b'} a,
[b,b']).$$ Cela d\Ž finit sur $W$ une structure d'alg\ bre
de Lie qui n'est autre que le \bf produit semi-direct \rm de
l'alg\ bre de Lie $B$ par l'alg\ bre de Lie $A[[Y]]$
relativement 
\ˆ \ l'homomorphisme $\sigma$ d\Ž fini ci-dessus (voir
Bourbaki [3], p. 17-18).

\

Bien entendu, l'alg\ bre de Lie $W$ ainsi construite ne d\Ž
pend que de $A$, de $B$ et de l'action $\rd : B \to S(Y)$,
mais pas de l'action $\rD : A \to S(X)$. On la d\Ž signera
par
$W(A,B;\rd)$ et on l'appellera \bf produit d'entrelacement
\rm de l'alg\ bre de Lie $B$ par l'alg\ bre de Lie $A$
relativement \ˆ \ l'action $\rd$.

\

\su{9 Action triangulaire}

\

Consid\Ž rons \ˆ \ nouveau une action $\rD$ de $A$ sur $X$,
une action $\rd$ de $B$ sur $Y$, le produit d'enrelacement
$W = W(A,B;\rd)$ et l'espace vectoriel produit $Z = X \times
Y$.

\

L'alg\ bre de Lie $S(Y)$ s'identifie naturellement \ˆ \ une
sous-alg\ bre de Lie de $S(Z)$ (lemme $<10>$).

\

On pose $T = S(X)$ et on consid\ re $T[[Y]]$ l'espace des
s\Ž ries formelles \ˆ \ variables dans $Y$ et coefficients
dans $T$. On l'identifie canoniquement \ˆ \ un sous-espace
de $S(Z)$ (th\Ž or\ me $<11>$).

\

Or, \ˆ \ chaque $a \in A[[Y]]$ et chaque $y \in Y$
correspond une s\Ž rie formelle $\rD_a \in T[[Y]]$ que l'on
identifie \ˆ \ l'\Ž l\Ž ment correspondant de $S(Z)$. Enfin
pour $(a,b) \in A[[Y]] \times B$, on pose
$$\Delta_{(a,b)} = \rD_a + \rd_b \ \ \text{un \Ž l\Ž ment de}
\ \ S(Z).$$

On montre (th\Ž or\ me $<12>$) que l'application ainsi d\Ž
finie
$$\Delta : W \to S(Z)$$
est une action formelle de $W$ sur $Z$. On
l'appellera \bf action triangulaire\rm, produit de l'action
$\rd$ par l'action $\rD$.

\

Pour 
$$a \in A[[Y]] \ , \ b \in B \ , \ x \in X \ , \ y \in
Y,$$
on donne un sens \ˆ \ l'\Ž galit\Ž 
$$\Delta_{(a,b)}(x,y) = \rD_{a(y)} + \rd_b(y).$$
On dira, de mani\ re imag\Ž e, que l'action de $W$ au point
$(x,y)$ est le r\Ž sultat de l'action de $B$ au point $y$ et
d'une action \dots \bf qui d\Ž pend de $y$ \rm \dots de $A$
au point $x$.

\

\su{10 Action fondamentale d'une alg\ bre de Lie sur
elle-m\ me}

\

Pour le crochet $[\xi,\eta] = \xi\eta - \eta\xi$ (voir
ci-dessus, au \bf 6\rm), on a vu que $S(B)$ est une alg\
bre de Lie.

\

On va d\Ž finir un homomorphisme canonique d'alg\ bres de
Lie
$$\rd : B \to S(B)$$
de la mani\ re suivante.

\

On commence par consid\Ž rer la \sl s\Ž rie g\Ž n\Ž ratrice
\rm
$$G(T) = \frac{T e^T}{e^T - 1} = \sum_{n \geq 0}
t_nT^n.\tag 1$$
Les coefficients $t_n$ appartiennent au corps $K$. Plus pr\Ž
cis\Ž ment,
$$t_0 = 1 \ , \ t_1 = 1/2 \ \ \text{et, pour} \ \ n \geq 1 \
, \ t_{2n} = \frac{b_{2n}}{(2n)!} \ , \ t_{2n+1} = 0,$$
o\ \ les $b_{2n}$ sont les nombres de \smc BERNOULLI\rm.

\

Pour $b \in B \ , \ y \in Y \ , \ s \in \N \ ,$ posons
$$\rd_{b,n}(y) = t_n (\rad \ y)^n(b)$$
o\ \ $\rad \ y : B \to B$ d\Ž signe l'application lin\Ž aire
adjointe
$$(\rad \ y)(b) = [y,b].$$
Ainsi, $\rd_{b,n}$ est un polyn\™ me homog\ ne de degr\Ž \
$n$ \ˆ \ variables et coefficients dans $B$ (lemme $<13>$).

\

On d\Ž signe par $\rd_b = (\rd_{b,n})$ la s\Ž rie formelle
correspondante. On d\Ž finit ainsi une application canonique
$\rd : B \to S(B)$, $b \mapsto \rd_b$, que l'on appellera
l'action \bf fondamentale \rm de $B$.

\

\su{Th\Ž or\ me} $<14>$. \sl Pour toute alg\ bre de Lie
$B$, l'action fondamentale $\rd : B \to S(B)$ est un
homomorphisme d'alg\ bres de Lie, autrement dit, $\rd$ est
une action formelle (\ˆ \ droite) de $B$ sur elle-m\ me.\rm

\

\su{11 Remarque}

\

Lorsque $B$ est une alg\ bre de Lie normable compl\ te sur
$\R$ ou $\C$, on sait lui associer \bf le groupuscule de Lie
d\Ž fini par $B$ \rm (voir Bourbaki [4], p. 168-169, ou
Kirillov [5]). Soit $G$ ce groupuscule ($G$ est un voisinage
de $0$ dans $B$). Bien entendu, l'alg\ bre de Lie $L(G)$ de
$G$ s'identifie \ˆ \ $B$.

\

Il existe, par hypoth\ se, \bf un morceau \rm de loi d'op\Ž
ration \ˆ \ doite analytique canonique du groupuscule $G$
sur la vari\Ž t\Ž \ $G$. Et, \ˆ \ ce \bf morceau \rm de
loi, correspond une \bf loi d'op\Ž ration infinit\Ž simale
\ˆ \ droite \rm analytique $\rD$ de $B = L(G)$ dans $G$
(voir Bourbaki [4], p. 165).

\

Dans ce cas, on montre (th\Ž or\ me $<15>$) que \bf l'action
fondamentale \rm de $B$ n'est autre que l'action \ˆ \ droite
de $B$ sur $B$ \bf d\Ž duite \rm de $\rD$ (voir ci-dessus,
l'exemple originel).

\

Pour $b \in B \ , \ y \in Y$, on a
$$\rd_b(y) = G(\rad \ y)(b) = \sum_{n \geq 0} t_n (\rad \
y)^n(b).\tag 2$$
On donne un sens \ˆ \ cette \Ž galit\Ž \ dans le cas g\Ž n\Ž
ral d'une alg\ bre de Lie quelconque.

\

\su{12 Produit d'entrelacement de deux alg\ bres de Lie}

\

On appellera \bf produit d'entrelacement \rm de l'alg\ bre
de Lie $B$ par l'alg\ bre de Lie $A$ le produit
d'entrelacement $W(A,B) = W(A,B;\rd)$ o\ \ $\rd$ est
\bf l'action fondamentale \rm de $B$.

\ 

\su{13 Repr\Ž sentation des extensions de $B$ par $A$ dans
le produit d'entrelacement $W(A,B)$} C'est l'analogue, pour
les alg\ bres de Lie, du premier th\Ž or\ me de
Kaloujnine-Krasner sur les groupes abstraits quelconques
(voir [6]).

\

Soit $A \to C \overset{\rp} \to \to B$ une extension de $B$
par $A$. Autrement dit, $\rp$ est un homomorphisme surjectif
de l'alg\ bre de Lie $C$ sur l'alg\ bre de Lie $B$, et $A$
est le noyau de $\rp$.

\

Soit $\rs : B \to C$ une application $K$-lin\Ž aire
quelconque telle que $\rp \circ \rs = \roman{id}_B$
(autrement dit, $\rs$ est une \bf section lin\Ž aire \rm de
$\rp$). On va associer \ˆ \ $\rs$ une application $f_{\rs} :
C \to W(A,B)$ que l'on appellera la \bf repr\Ž sentation
associ\Ž e \rm \ˆ \ $\rs$.

\

Pour $c \in C \ , \ y \in B \ , \ m \in \N$, on pose $z = \rs
(y)$ et 
$$h_{c,m} = \frac{1}{m!} (\rad \ z)^m(c) - \sum_{n+r =
m} \frac{t_r}{(n+1)!} (\rad \ z )^n(\rs \circ \rp) (\rad \
z)^r(c)\tag 3$$
o\ \ les coefficients $t_r$ sont d\Ž finis par la relation
(1) ci-dessus. On montre (lemme $<16>$) que $h_{c,m}(y)
\in A$. De sorte que $h_{c,m} \in A[B]_m$. Aisni
$h_c = (h_{c,m}) \in A[[B]]$.

\

On pose enfin $f_\rs(c) = (h_c,\rp(c))$.

\

\su{Th\Ž or\ me} $<17>$. \sl Soit $A \to C \overset{\rp} \to
\to B$ une extension de $B$ par $A$. Pour toute section
lin\Ž aire $\rs$ de $\rp$, la repr\Ž sentation associ\Ž e
$f_{\rs} : C \to W(A,B)$ est un homomorphisme injectif
de l'alg\ bre de Lie $C$ dans l'alg\ bre de Lie
$W(A,B)$.\rm

\

\su{Pr\Ž cisions} La premi\ re annonce de ces r\Ž sultats a
\Ž t\Ž
\ faite sous forme d'un rapport pr\Ž liminaire dans les
\lq\lq Abstracts of the A. M. S.", sous la r\Ž f\Ž rence
85T-27-237. Plus tard, une version de ce texte a \Ž t\Ž \
publi\Ž e dans
$E\Lambda E\Upsilon\Theta EPIA$ (ELEFTERIA) 3 (1985)
290-304.

\

\su{Nota bene} La notion de produit d'entrelacement d\Ž
finie ci-dessus diff\ re essentiellement de celle de \lq\lq
verbal $\frak V$-wreath product" introduite par A. L. \u
SMELK\'IN dans \sl Trans. Moscow Math. Soc. \rm \bf 29 \rm
(1973) p. 239-252.

\

\

\

\centerline{\bf Voici une version courte en anglais\rm}

\

\

\

\centerline{\bf Wreath products and triangular actions of
Lie algebras\rm}

\

\head Abstract\endhead

\

Formal actions of Lie algebras over
vector spaces are introduced in a purely algebraic way, as a
mimic of infinitesimal operations of Banach Lie algebras over
Banach analytic manifolds. In analogy with the case of
abstract groups, complete wreath products and triangular
actions are then defined for Lie algebras acting \sl en
cascade \rm over vector spaces. Finally, a
Kaloujnine-Krasner type theorem for Lie algebra extensions
is proved. 

\

\head A moderately detailed english summary\endhead

\

\centerline{\bf \lq\lq Can wreath products for Lie groups be
defined ?"}

\

That was the question that M. Krasner once put to the
first-named author.

\

The problem of defining the \sl true \rm wreath product of
two Lie groups is not easy, and seems to be still open.

\

In order to get around the obstacle, and as a first step, we
tried to find a good definition for the wreath product of two
Lie \bf algebras\rm. That is what the present note is about.

\

A paper giving the details is in preparation. The numbers
such as
$<n>$ refer to that coming paper.

\

\su{A brief sketch for the definition of the wreath product
for Lie algebras} 

\

All vector spaces and algebras are over a given field $K$ of
characteristic zero.  

\

Let $A$ and $B$ be
Lie algebras. Set
$$A[B]_n = \{ f : B \to A : \exists u : B^n \to A  \ \
\text{an $n$-linear map such that} \ \ f(x) =
u(x,\dots,x)\},$$
$$A[[B]] = \prod_{n > 0} A[B]_n \ , \ S(B) = B[[B]].$$
A given $f \in A[B]_n$ is said to be \bf determined \rm  by
$u : B^n \to A$ whenever
$f(x) = u(x,\dots,x)$. 

\

\noi For $f_p \in B[B]_p$, determined by $u_p$, and $g_q \in
B[B]_q$, define
$$(g_q.f_p)(x) = \sum_{1 \leq i \leq p}
u_p(x,\dots,g_q(x),\dots,x) \ , \ g_q(x) \ \ \text{in the
i-th place}.$$ 
For $f = (f_p)$ and $g = (g_q)$ in
$S(B)$, define 
$$[f,g] = (h_n) \ , \ h_n(x) =
\sum_{p+q = n+1}(f_p.g_q - g_q.f_p)(x).$$
For $f = (f_p)$ and $g = (g_p)$ in
$A[[B]]$, define
$$[f,g] =
(h_n)\ , \ h_n(x) = \sum_{p+q=n}[f_p(x),g_q(x)].$$ 

\

\su{For the brackets defined above, $A[[B]]$ and
$S(B)$ are Lie algebras, ($<7>$)} 
\

\noi Next, define the $t_n$'s by their \sl generating
function
\rm : 
$$G(T) = \sum t^n T^n = \frac{T e^T}{e^T - 1}.$$
Of course,
$$t_0 = 1 \ , \ t_1 = 1/2 \ \ \text{and, \ for} \ \ n \geq 1
\ , \ t_{2n} = \frac{b_{2n}}{(2n)!} \ , \ t_{2n+1} = 0,$$
where the $b_{2n}$'s are the \smc BERNOULLI \rm numbers.
Then define $\rd : B \to S(B)$, setting
$$\rd_b (y) = \sum t_n (\rad \ y)^n(b).$$
For any $b \in B$, the map $a \mapsto \rd_b.a = \sum
d_{b,q}.a_p$ is a derivation of the Lie algebra $A[[B]]$;
and the map $b \mapsto
\rd_b$ is a Lie algebra homomorphism $\sigma : B \to \frak
d(A[[B]]$ of
$B$ into the derivation algebra of $A[[B]]$, ($<14>$).

\

\su{Definition}  We define the \bf wreath product \rm \sl of
the Lie algebra $B$ by the Lie algebra 
$A$ \rm to be the semi-direct product $W(A,B) = A[[B]]
\underset
\sigma \to\times B$ relative to $\sigma$.

\

\su{Representation of  Lie algebra
extensions} 

\

Given any exact sequence
$0
\to A
\to C
\overset \rp \to\to B \to 0$ of Lie algebras, any linear
section $\rs : B \to C$ of $\rp$, (i.e. $\rp \circ \rs =
\roman{id}_B$), any elements
$c
\in C$ and
$y
\in B$, set $z = \rs(y)$ and 
$$h_{c,m} = \frac{1}{m!} (\rad \ z)^m (c) -
\sum_{n+r = m} \frac{t_r}{(n+1)!} (\rad \ z )^n (\rs \circ
\rp) (\rad \ z)^r(c),$$
$$h_c = (h_{c,m}) \ , \ f(c) = (h_c,\rp(c)).$$
\bf Then $f : C \to W(A,B)$ is an injective Lie algebra
homomorphism\rm, ($<17>$).

\

Moreover, general actions of Lie algebras over
vector spaces are introduced, as a new notion. We also
exhibit a natural \sl triangular \rm action of the wreath
product $W(A,B)$ over the vector product space $A
\times B$.

\

Here are a few hints about those generalizations.

\

\head More details\endhead

\

All vector spaces and algebras are over a given field $K$ of
characteristic zero. In the sequel, $E,F,X,Y,$ are vector
spaces,
$A,B,C,$ are Lie algebras, while $m,n,r,$ are integers. 

\

In order to define the wreath product of two Lie algebras, we
needed a notion of an
\bf action
\rm for Lie algebras which had to be general enough for our
needs. We were thus lead to introduce a certain number of
new tools (which might also be interesting for their own
sake). Here
are some of them.

\

We introduced, namely, the notion of a derivation relative
to \sl a formal series with variables and coefficients in a
vector space $X$\rm. This leads to the introduction
of a new object \bf the Lie algebra $S(X)$ of formal seies
on $X$\rm.

\

For each $m$, $L_m(E;F)$ is the set of $m$-linear maps
$E \to F$.

\

\su{1 Formal series} A map $f : E
\to F$ is called a  \bf homogeneous polynomial having degree
$m$
with
variables in $E$ and coefficients in $F$ \rm whenever there
is a $u \in L_m(E;F)$ such that
$$f(x) = u(x,x,\dots,x) \ \ \text{for each} \ \ x \in E.$$
Then, $f$ is said to be \bf determined \rm by $u$.

\

Let $F[E]_m$ be the set of all those homogeneous polynomials
having degree \ $m$. This set is, naturally, endowed with a
structure of vector space (over
$K$). Let
$$F[E] = \underset m\to \bigoplus F[E]_m \ , \ F[[E]] =
\prod_m F[E]_m.$$
We then define a \bf polynomial \rm (resp. a \bf formal
series\rm)
with variables in $E$ and coefficients in $F$ to be any 
element in $F[E]$ (resp. in F[[E]]).

\

\su{2 Symmetrisation} For each
$$u \in L_m(E;F) \ , \ z = (z_1,\dots,z_r) \in E^r \ , \ p =
(p_1,\dots,p_r) \in \Z^r,$$
if $m = p_1 + \dots + p_r$, let $\tilde
u(z;p)$ be the sum of all terms having the form
$u(x_1,\dots,x_m)$ where, among the $x_i$'s $(i =
1,\dots,m)$, exactly $p_j$ are equal to $z_j$. If no such
terms exist, just set $\tilde
u(z;p) = 0$.

\

\su{3 When $u$ and $v$ both determine the same homogeneous
polynomial, then $\tilde u(z;p) = \tilde v(z;p) \ \
\text{for all} \ \ z
\in E^r \ , \ p \in \Z^r$, ($<1>$)}

\

\su{5 Derivation relative to a formal series}

\

Let $\xi \in X[X]_r$ and $f \in F[X]_m$.

\ 

Suppose that $f$ is \bf determined \rm  by $u \in
L_m[X;F)$. For each $x \in X$, let $(\xi
f)(x)$ be the sum of all terms of the form
$u(x_1,\dots,x_m)$ where, among the $x_i$'s $(i =
1,\dots,m)$, only one is equal to $\xi(x)$ and the others
are all equal to
$x$. That is,
$$(\xi f)(x) = \tilde
u((\xi(x),x);(1,m-1)).$$ 
Whence a map $\xi f : X \to F$.

\

We show that  $\xi f$, a homogeneous polynomial ($<3>$) which
does not depend on the choice of a particular
$u$, is such that $\xi f \in
F[X]_s$ where $s = r + m - 1$, and it is agreed that
$F[X]_s =
\{0\}$ for $s < 0$.

\

\su{Definition} Write $S(X)$ for $X[[X]]$.

\

Let two formal series $\xi =
(\xi_r)
\in S(X)$ and
$f = (f_m) \in F[[X]]$ be given. For each $s \geq
0$, let
$$g_s = \sum_{r + m - 1 = s} \xi_r f_m.$$
Denote by $\xi f$ the formal series $(g_s) \in
F[[x]]$, and call it the \bf derivative of
$f$ \sl relative to the formal series
$\xi$\rm.

\

\su{6 The Lie algebra $S(X)$}  Given formal series
$\xi$ and $\eta$ in
$S(X)$, take $\xi \eta$, the derivative of 
$\eta$ relative to $\xi$, and $\eta \xi$, the derivative of
$\xi$ relative to $\eta$. Set $[\xi,\eta] = \xi \eta -
\eta
\xi$.

\

\su{The vector space 
$S(X)$ with the bracket thus defined on it is a Lie
alebra, ($<4>$)} 
\

\su{7 Action of a Lie algbra over a vector space}
A \bf formal action \rm of a Lie
algebra $A$ over a vector space $X$ is defined to be any
Lie algebra homomorphism $\rD : A \to S(X)$. Thus, the Lie
algebra $A$ acts over the vector space $X$ \sl through \rm
its Lie algebra $S(X)$.

\

\su{The original example} Take $K = \R$ or $\C$.
Suppose
$A$ is a complete normable algebra and $X$
a Banach space. Take an open neighbourhood $U$ of
$0$ in $X$, and an analytic \bf infinitesimal operation
law\rm, $a \mapsto \rD_a$, of $A$ in the analytic
manifold $U$ (see Bourbaki [4], p.139). 

\

Thus, for each $a \in A$, the vector field $\rD_a$
on
$U$ is analytic. Of course, the natural embedding $h :
U \to X$ is analytic. Consider the map
$\rD_a(h) : U \to X$ (see Bourbaki [2], 8.2.2 et 8.2.3,
p.10).

\

Since it is analytic,
it is represented as a (convergent) formal series, having
continuous components, in the neighbourhood of $0$, that is,
an element of $S(X)$, which we still denote
$\rD_a$.

\

\su{The map $\rD : A
\to S(X)$ thus defined is a formal action of $A$
over $X$, ($<5>$)}
 
\

\su{8 Wreath products} Let a formal action
$\rD$ of $A$ over $X$ and a formal action
$\rd$ of $B$ over $Y$ be given. Take the two vector spaces
$$W = A[[Y]] \times B \ \ \text{and} \ \  \ Z = X \times Y.$$

\

We first endow $W$ with a Lie algebra structure, which we
call the \bf wreath product\rm,
and then define a formal action $\Delta : W \to S(Z)$, which
we call the \bf triangular action\rm.

\

\su{8.1 The bracket on $A[[Y]]$} 

\

For $f \in A[Y]_n$ and $g \in
A[Y]_r$ and each $y \in Y$, set $[f,g](y) =
[f(y),g(y)]$. This defines a map $[f,g] : Y \to
A$, and  $[f,g] \in A[Y]_{n + r}$, ($<6>$).
More generally, for $f = (f_n) \in A[[Y]]$ and $g =
(g_r) \in A[[Y]]$, set $[f,g]_s = \sum_{n+r =s}
[f_n,g_r]$ and also $[f,g] =([f,g]_s) \in A[[Y]]$. This is a
Lie algebra bracket on $A[[Y]]$, ($<7>$).

\

Let $\frak d (A[[Y]])$ be the Lie algebra of the derivations
of the Lie algebra $A[[Y]]$.

\

\su{8.2 A homomorphism from $B$ into $\frak d (A[[Y]])$}

\

Take the formal action $\rd : B \to S(Y)$.
For each $b \in B$, we have $\rd_b \in S(Y)$; and for each
$a \in A[[Y]]$, the derivative $\rd_b a$ of $a$ relative to
$\rd_b$ belongs to $A[[Y]]$. So $\rd_b$ is a map
from $A[[Y]]$ into itself. This map is a derivation of the
Lie algebra
$A[[Y]]$ which we denote as $\sigma(b)$, ($<8>$). This map
$\sigma : B \to \frak d (A[[Y]])$ is a Lie algebra
homomorphism, ($<9>$). 
\

\su{8.3 The bracket on $W = A[[Y]] \times B$}

\

Let $(a,b)$ et $(a',b')$ be elements of $W$.

\

Set
$$[(a,b),(a',b')] = ([a,a'] + \rd_b a' - \rd_{b'} a,
[b,b']).$$ 
This turns $W$ into a Lie algebra which is the \bf 
semi-direct product \rm of the Lie algebra $B$ by the Lie
algebra
$A[[Y]]$ relative to the homomorphism $\sigma$, (see
Bourbaki [3], p. 17-18).

\

Of course, this Lie algebra $W$ depends only
on $A$, $B$, and the action $\rd : B \to S(Y)$,
but not on  $\rD : A \to S(X)$. We denote it as
$W(A,B;\rd)$ and call it \bf the wreath product 
\rm of the Lie algebra $B$ by the Lie algebra $A$
relative to the action $\rd$.

\

\su{9 The triangular action}

\

Take again an action $\rD$ of $A$ over $X$,
an action $\rd$ of $B$ over $Y$, the wreath product
$W = W(A,B;\rd)$ and the vectorial product space $Z = X
\times Y$.

\

The Lie algebra $S(Y)$ is identified, in a natural way, to a
Lie subalgebra of $S(Z)$, ($<10>$). 

\

Let $T = S(X)$ and take $T[[Y]]$, the space of formal series
with variables in $Y$ and coefficients in $T$,
canonically identified to a subspace of $S(Z)$, ($<11>$). 

\

To each $a \in A[[Y]]$ and each $y \in Y$ there
corresponds a formal series $D_a \in T[[Y]]$ identified to
the corresponding element in $S(Z)$. Now, for $(a,b) \in
A[[Y]] \times B$, let
$$\Delta_{(a,b)} = \rD_a + \rd_b \ \ \text{an element in}
\ \ S(Z).$$

The map $\Delta : W \to S(Z)$
is a formal action of $W$ over $Z$,  the \bf 
triangular action\rm, a product of the actions
$\rd$ and $\rD$, ($<12>$). 

\

For 
$$a \in A[[Y]] \ , \ b \in B \ , \ x \in X \ , \ y \in
Y,$$
the following equality makes sense :
$$\Delta_{(a,b)}(x,y) = \rD_{a(y)} + \rd_b(y).$$
One can say, figuratively, that the action of $W$
at point
$(x,y)$ is the result of the action of $B$ at point $y$ and
an action \dots \bf which depends on $y$ \rm \dots of
$A$ at point $x$.

\

\

\

\

\

\

\


\Refs\nofrills{Bibliographie}

\ref \no1 \by Bourbaki, N. \book Vari\Ž t\Ž s diff\Ž
rentielles et analytiques, Fascicule des r\Ž sultats,
paragraphes 1 \ˆ \ 7
\publ Hermann, Paris \yr 1967 \endref

\

\ref \no2 \bysame \book  Vari\Ž t\Ž s diff\Ž
rentielles et analytiques, Fascicule des r\Ž sultats,
paragraphes 8 \ˆ \ 15 \publ Hermann, Paris \yr 1971 \endref

\

\ref \no3 \bysame \book Groupes et alg\ bres de Lie,
Chapitre 1 \publ Hermann, Paris \yr 1972 \endref

\

\ref \no4 \bysame \book Groupes et alg\ bres de Lie,
Chapitres 2 et 3 \publ Hermann, Paris \yr 1972 \endref


\

\ref \no5 \by Kirillov, A. \book El\Ž ments de la th\Ž orie
des repr\Ž sentations, traduction fran\c caises \publ
Editions Mir , Moscou \yr 1974 \endref

\

\ref \no6 \by Neumann, H. \inbook Varieties of groups
(theorem 22.21) \pages 45-46
\publ Springer , Berlin \yr 1967 \endref

\

\endRefs
\enddocument